\documentclass[11pt]{article}
\usepackage{amscd}
\usepackage{amsfonts}
\usepackage{amsmath}
\usepackage{amssymb}
\usepackage{amsthm}
\usepackage{bbm}
\usepackage{CJK}
\usepackage{fancyhdr}
\usepackage{graphicx}
\usepackage{hyperref}
\usepackage{indentfirst}
\usepackage{latexsym}
\usepackage{mathrsfs}
\usepackage{xypic}

\usepackage[top=1in,bottom=1in,left=1.25in,right=1.25in]{geometry}
\theoremstyle{definition}

\def\<{\langle}
\def\>{\rangle}
\def\a{\alpha}

\def\c{\cdot}

\def\D{\Delta}

\def\o{\otimes}

\def\v{\varepsilon}

\date{}
\begin{document}
\renewcommand{\baselinestretch}{1.2}
\renewcommand{\arraystretch}{1.0}
\title{\bf Yetter-Drinfeld modules for weak Hom-Hopf algebras
}
\date{}
\author {{\bf Shuangjian Guo$^{1}$,  Yizheng Li$^{1}$,  Shengxiang Wang$^{2}$}\\
{\small 1. School of Mathematics and Statistics, Guizhou University of Finance and Economics}\\
{\small Guiyang,  550025, P. R. China}\\
{\small Email: shuangjianguo@yahoo.com,  Email: yzli1001@163.com }\\
{\small 2. School of Mathematics and Finance, Chuzhou University}\\
{\small Chuzhou 239000,   P. R. China}\\
{\small E-mail: wangsx-math@163.com}}
 \maketitle
\begin{center}
\begin{minipage}{12.cm}
\begin{center}{\bf ABSTRACT}\end{center}
  The aim of this paper is to define and study  Yetter-Drinfeld modules over  weak Hom-Hopf algebras.
We  show that the category ${}_H{\cal
WYD}^H$ of Yetter-Drinfeld modules  with bijective structure
maps over  weak Hom-Hopf algebras is a rigid category  and a braided monoidal category,  and obtain a new solution of quantum Hom-Yang-Baxter equation. It turns out that,  If $H$ is quasitriangular (respectively,
coquasitriangular)weak Hom-Hopf algebras,   the category of
modules (respectively, comodules) with bijective structure maps over $H$ is a braided monoidal subcategory of the category ${}_H{\cal
WYD}^H$ of  Yetter-Drinfeld modules over  weak Hom-Hopf algebras.

 \vskip 0.5cm

\end{minipage}
\end{center}
\section*{0. Introduction}

The first examples of Hom-type algebras were related to $q$-deformations of Witt and Virasoro
algebras, which play an important role in Physics, mainly in conformal field theory. The $q$-deformations of Witt
and Virasoro algebras are obtained when the derivation is replaced by a $\sigma$-derivation. It was observed
in the pioneering works (See \cite{CE90}-\cite{H99}).  Motivated by these examples and their generalization,   Hartwig, Larsson and Silvestrov introduced the Hom-Lie algebras when they concerned about the $q$-deformations of Witt and Virasoro algebras
in \cite{JDS}.   In a Hom-Lie algebra, the Jacobi identity is replaced by the so called Hom-Jacobi identity via an homomorphism. Hom-associative algebras,
the corresponding structure of associative algebras, were introduced by Makhlouf and Silvestrov in \cite{MS08}. The associativity of the Hom-algebra
is twisted by an endomorphism (here we call it the Hom structure map). The
generalized notions, Hom-bialgebras, Hom-Hopf algebras were developed in \cite{MS09}, \cite{MS10}. Caenepeel and Goyvaerts studied in \cite{CG11} Hom-bialgebras and Hom-Hopf algebras from
a categorical view point, and called them monoidal Hom-bialgebras and monoidal Hom-Hopf algebras respectively, which are slightly different from the above Hom-bialgebras
and Hom-Hopf algebras.  Thus a monoidal Hom-bialgebra is Hom-bialgebra if and only if the Hom-structure map $\a$ satisfies
$\a^2 = id$.  Yau introduced Quasitriangular Hom-bialgebras
in \cite{Y09}), which provided a solution of the quantum Hom-Yang-Baxter euqation,  a twisted version of the quantum Yang-Baxter equation called
the Hom-Yang-Baxter equation  in \cite{Y11}.     Zhang and Wang
introduced weak Hom-Hopf algebra $H$,  which is generalization of both Hom-Hopf algebras and weak Hopf algebras, and discussed
the  category $Rep(H)$ (resp. $Corep(H)$) of Hom-modules (resp. Hom-comodules) with bijective Hom-structure maps,   they  proved that if $H$ is a (co)quasitrialgular weak Hom-bialgebra (resp. ribbon weak Hom-Hopf algebra),
then $Rep(H)$ (resp. $Corep(H)$) is a braided monoidal category (resp. ribbon category) in \cite{ZW}.

Makhlouf and Panaite defined and studied Yetter-Drinfeld modules over Hom bialgebras,
a generalized version of bialgebras obtained by modifying the algebra and coalgebra structures by a homomorphism. Yetter-Drinfeld modules over a Hom bialgebra with bijective
structure map provide solutions of the Hom-Yang-Baxter equation in \cite{MP14}. It is well known that the Yetter-Drinfeld modules category of
a (weak) Hopf algebra is a rigid monoidal category, and is braided. Does this result remain true in a weak Hom-Hopf algebra? How the corresponding results appear
under the condition that the associativity and coassociativity are twisted by an endomorphism? Is there any relation between this Yetter-Drinfeld modules category and module category or
comodule category? This is the motivation of the present article. In order to investigate these questions, we introduce the definition of  Yetter-Drinfeld modules over  weak Hom-Hopf algebras,  which is generalization of both weak
Yetter-Drinfeld modules introduced by \cite{CW05} or \cite{N02} and Hom-Yetter-Drinfeld modules introduced by \cite{MP14} or \cite{MP15}, and consider that when the  Yetter-Drinfeld modules category of a weak Hom-Hopf algebra is is a rigid monoidal category, and is braided.

To make sure that the  Yetter-Drinfeld modules category of a weak Hom-Hopf algebra $H$ is a monoidal category, we need $H$ is unital and counital, and the Hom structure maps over $H$  are all bijective maps.

The paper is organized as follows. In Section 2,  we recall now several concepts and results, fixing thus the terminology to be used in the rest of
the paper.

 In Section 3, we  introduce the definition of  Yetter-Drinfeld modules over  weak Hom-Hopf algebras and
 show the category ${}_H{\cal
WYD}^H$ of  Yetter-Drinfeld modules is a monoidal category and a rigid category.

 In Section 4, we  show that the category ${}_H{\cal
WYD}^H$ of  Yetter-Drinfeld modules is a braided monoidal category and  obtain a new solution of quantum Hom-Yang-Baxter equation. It turns out that, if $H$ is a quasitriangular weak Hom-Hopf algebra, the category
of left $H$-modules with bijective structure maps is a braided monoidal subcategory of the category ${}_H{\cal
WYD}^H$ of Yetter-Drinfeld modules.

 In Section 5, we find another braided monoidal category structure on the category ${}_H{\cal
WYD}^H$ of  Yetter-Drinfeld modules , with the property
that if $H$ is a coquasitriangular weak Hom-Hopf algebra, then ${}_H{\cal
WYD}^H$ contains the category of right $H$-comodules
with bijective structure maps as a braided monoidal category.

\section*{2. Preliminaries}
\def\theequation{2. \arabic{equation}}
\setcounter{equation} {0} \hskip\parindent

Throughout the paper, we let $\Bbbk$ be a fixed
 field and all algebras are supposed to be over $\Bbbk$. For the comultiplication
 $\D $ of a vector space $C$, we use the Sweedler-Heyneman's notation:
 $$
 \Delta(c)=c_{1}\o c_{2},
 $$
for any $c\in C$. $\tau$ means the flip map $\tau(a \o b) = b \o a$.
When we say a "Hom-algebra" or a "Hom-coalgebra", we mean the unital Hom-algebra and counital Hom-coalgebra.

In this section, we will review several definitions and notations related to weak Hom-Hopf algebras and rigid categories.

\vskip 0.5cm
 {\bf 2.1. Hom-algebras and Hom-coalgebras.}
\vskip 0.5cm

Recall from \cite{MS08} that
a \emph{Hom-associative algebra} is a quadruple $(A,\mu,\eta,\alpha_A)$, in which $A$ is a linear space, $\alpha_A:A\rightarrow A$, $\mu: A \o A\rightarrow A$ and $\eta:\Bbbk \rightarrow A$ are linear maps, with notation $\mu(a \o b) = ab$ and $\eta(1_\Bbbk) = 1_A$, satisfying the following conditions, for all $a,b,c \in A$:\\
$\left\{\begin{array}{l}
(1)~~\alpha_A(ab) = \alpha_A(a)\a(b);\\
(2)~~\alpha_A(a)(bc) = (ab)\alpha_A(c);\\
(3)~~\alpha_A(1_A) = 1_A;\\
(4)~~1_Aa = a1_A = \alpha_A(a).
\end{array}\right.$
\\

\emph{A morphism $f:A\rightarrow B$ of Hom-algebras} is a linear map such that $\a_B\circ f = f\circ \a_A$, $f(1_A) = 1_B$ and $\mu_B\circ (f \o f) = f\circ \mu_A$.

Let $A$ be a Hom-algebra. Recall that a \emph{left $A$-module} is a triple $(M,\a_M, \theta_M)$, where $M$ is a $\Bbbk$-space, $\a_M:M\rightarrow M$ and $\theta_M:A \o M\rightarrow M$ are linear maps with notation $\theta_M(a \o m) = a\cdot m$, satisfying the following conditions, for all $a,b\in A, m\in M$:\\
$\left\{\begin{array}{l}
(1)~~\a(a\cdot m) = \a(a)\cdot\a_M(m);\\
(2)~~\a(a)\cdot(b\cdot m) = (ab)\cdot\a_M(m);\\
(3)~~1_A\cdot m = \a_M(m).
\end{array}\right.$
\\

\emph{A morphism $f:M\rightarrow N$ of $A$-modules} is a linear map such that $\a_N\circ f = f\circ \a_M$ and $\theta_N\circ (id_A \o f) = f\circ \theta_M$.

Recall from \cite{MS10} that
a \emph{Hom-coassociative coalgebra} is a quadruple $(C,\Delta,\varepsilon,\a_C)$, in which $C$ is a linear space, $\a:C\rightarrow C$, $\Delta: C\rightarrow C\o C$ and $\varepsilon:C\rightarrow \Bbbk$ are linear maps, with notation $\Delta(c) = c_1 \o c_2$, satisfying the following conditions for all $c \in C$:\\
$\left\{\begin{array}{l}
(1)~~\Delta(\a_C(c)) = \a_C(c_1)\a(c_2);\\
(2)~~\a_C(c_1) \o \Delta(c_{2}) = \Delta(c_1) \o \a_C(c_2);\\
(3)~~\varepsilon \circ \a_C = \varepsilon;\\
(4)~~\varepsilon(c_1)c_2 = c_1\varepsilon(c_2) = \a_C(c).
\end{array}\right.$
\\

\emph{A morphism $f:C\rightarrow D$ of Hom-coalgebras} is a linear map such that $\a_D\circ f = f\circ \a_C$, $\varepsilon_C = \varepsilon_D\circ f$ and $\Delta_D\circ f = (f \o f)\circ \Delta_C$.

Let $C$ be a Hom-coalgebra. Recall that a \emph{right $C$-comodule} is a triple $(M,\a_M, \rho_M)$, where $M$ is a $\Bbbk$-space, $\a_M:M\rightarrow M$ and $\rho_M:M\rightarrow M \o C$ are linear maps with notation $\rho_M(m) = m_{(0)}\o m_{(1)}$, satisfying the following conditions for all $m \in M$:\\
$\left\{\begin{array}{l}
(1)~~\rho_M(\a_M(m)) = \a_M(m_{(0)})\o \a_C(m_{(1)});\\
(2)~~\a_M(m_{(0)})\o \Delta(m_{(1)}) = \rho_M(m_{(0)}) \o \a_C(m_{(1)});\\
(3)~~\varepsilon(m_1)m_0 = \a_M(m).
\end{array}\right.$
\\

\emph{A morphism $f:M\rightarrow N$ of $C$-comodules} is a linear map such that $\a_N\circ f = f\circ \a_M$ and $\rho_N\circ f=(id_C \o f)\circ \rho_M$.

 Recall from \cite{ZW} that a \emph{weak Hom-bialgebra} is a sextuple  $H=(H,\a_H,\mu,\eta,\Delta,\varepsilon)$   if $(H,\a_H)$ is both a Hom-algebra and a Hom-coalgebra, satisfying the following identities for any $a,b,c \in H$:\\
$\left\{\begin{array}{l}
(1)~~\Delta(ab)  = \Delta(a)\Delta(b);\\
(2)~~\varepsilon((ab)c) = \varepsilon(ab_1)\varepsilon(b_2c),~~\varepsilon(a(bc)) = \varepsilon(ab_2)\varepsilon(b_1c);\\
(3)~~(\Delta \o id_H)\Delta(1_H) = 1_1 \o 1_21'_{1} \o 1'_{2},~~(id_H \o \Delta)\Delta(1_H) = 1_1 \o1'_{1}1_2 \o 1'_{2}.
\end{array}\right.$

Recall from \cite{ZW} that
a \emph{Weak Hom-Hopf algebra} is a septuple $(H,\mu$, $\eta$, $\Delta$, $\varepsilon,S,\a_H)$, in which $(H, \a_H)$ is a weak Hom-bialgebra,  if $H$ endowed with a $\Bbbk$-linear map $S$ (the \textsl{antipode}), such that for any $h,g \in H$, the following conditions hold:
\\
$\left\{\begin{array}{l}
(1)~~S\circ \a_H = \a_H\circ S;\\
(2)~~h_1S(h_2) = \varepsilon_t(h),~~S(h_1)h_2 = \varepsilon_s(h);\\
(3)~~S(hg) = S(g)S(h),~~S(1_H) = 1_H;\\
(4)~~\Delta(S(h)) = S(h_2) \o S(h_1),~~\varepsilon\circ S = \varepsilon.
\end{array}\right.$

Let $(H,\a_H)$ be a weak Hom-bialgebra. Define linear maps $\varepsilon_t$ and $\varepsilon_s$ by the formula
$$
\varepsilon_t(h) = \varepsilon(1_{1}h)1_{2},~~~\varepsilon_s(h) = 1_{1}\varepsilon(h1_{2}),
$$
for any $h \in H$, where $\varepsilon_t$, $\varepsilon_s$ are called the \textsl{target} and \textsl{source counital maps}. We adopt the notations $H_t = \varepsilon_t(H)$ and $H_s = \varepsilon_s(H)$ for their images.

Similarly, we define the linear maps $\widehat{\varepsilon_t}$ and $\widehat{\varepsilon_s} $ by the formula
$$
\widehat{\varepsilon_t}(h) = \varepsilon(h 1_{1})1_{2},~~~\widehat{\varepsilon_s}(h) = 1_{1}\varepsilon(1_{2}h),
$$
for any $h \in H$. Their images are denoted by $\widehat{H_t} = \widehat{\varepsilon_t}(H)$ and $\widehat{H_s} = \widehat{\varepsilon_s}(H)$.

\vskip 0.5cm
 {\bf 2.2. Duality and rigid categories.}
\vskip 0.5cm

Recall from \cite{K95} that let $(\mathcal {C}, \o,I, a,l,r)$ be a monoidal category. $V \in \mathcal {C}$, a \emph{left dual} of $V$ is a triple $(V^\ast, ev_V, coev_V)$, where $V^\ast$ is an object, $ev_V : V^\ast \o V\rightarrow I$ and $coev_V : I\rightarrow V \o V^\ast$ are morphisms in $\mathcal {C}$, satisfying
$$ r_V\circ(id_V \o ev_V)\circ a_{V,V^\ast,V} \circ (coev_V \o id_V)\circ l_V^{-1} = id_V,$$
and
$$l_{V^\ast}\circ (ev_V \o id_{V^\ast})\circ a^{-1}_{V^\ast,V,V^\ast}\circ (id_{V^\ast} \o coev_V)\circ r_{V^\ast}^{-1} = id_{V^\ast}.$$

Similarly, a \emph{right dual} of $V$ is a triple $({}^\ast V, \widetilde{ev}_V, \widetilde{coev}_V)$, where ${}^\ast V$ is an object, $\widetilde{ev}_V : V \o {}^\ast V\rightarrow I$ and $\widetilde{coev}_V : I\rightarrow {}^\ast V \o V$ are morphisms in $\mathcal {C}$, satisfying
$$ r_{{}^\ast V}\circ(id_V \o \widetilde{ev}_V)\circ a_{{}^\ast V,V,{}^\ast V} \circ (\widetilde{coev}_V \o id_{{}^\ast V})\circ l_{{}^\ast V}^{-1} = id_{{}^\ast V},$$
and
$$l_{V}\circ (\widetilde{ev}_V \o id_{V})\circ a^{-1}_{V,{}^\ast V,V}\circ (id_{V} \o \widetilde{coev}_V)\circ r_{V}^{-1} = id_{V}.$$

If each object in $\mathcal {C}$ admits a left dual (resp. a right dual, both a left dual and a right dual), then $\mathcal {C}$ is called a \emph{left rigid category} (resp. a \emph{right rigid category}, a \emph{rigid category}).

\section*{3. Left-right Yetter-Drinfeld modules over a weak Hom-Hopf algebra}
\def\theequation{3. \arabic{equation}}
\setcounter{equation}{0} \hskip\parindent

{\bf Definition 3.1.} Let $(H,\a_H)$ be a weak Hom-Hopf algebra.
 A Yetter-Drinfeld module over $H$ is a vector space $(M,\a_M)$,
such that $M$ is a unital left $H$-module (with notation $h\o m\mapsto
h\cdot m$) and a counital right $H$-comodule (with notation $m\o h \mapsto
m_{(0)}\o m_{(1)}$) with the following compatibility condition:
\begin{equation}
(h\c m)_{(0)}\o (h\c m)_{(1)}=\a^{-1}_H(h_{21})\c m_{(0)}\o
[\a^{-2}(h_{22})\a_H^{-1}(m_{(1)})]S^{-1}(h_{1}),
\end{equation}
for all $h\in H$ and $m\in M$. We denote by ${}_H{\cal
WYD}^H$  the category of Yetter-Drinfeld
modules, morphisms being the $H$-linear $H$-colinear maps.

{\bf Proposition  3.2.}  one has that(3.1) is equivalent to the following equations
\begin{eqnarray}
\rho(m)= m_{(0)}\o m_{(1)}\in M\o_t H\triangleq (1_{1}\o 1_{2})\c ( M\o H)\\
\a_H(h_{1})\c m_{(0)}\o \a_H^{2}(h_{2})\a_H(m_{(1)})=(h_{2}\c m)_{(0)} \o (h_{2}\c
m)_{(1)}\a_H^{2}(h_{1})
\end{eqnarray}
{\it Proof }   $(3.1)\Longrightarrow(3.2), (3.3)$. We have
\begin{eqnarray*}
&&m_{(0)}\o m_{(1)}\\
   &=&\a_H^{-1}(1_{21})\c \a_M^{-1}(m_{(0)})\o[\a_H^{-2}(1_{22})\a_H^{-2}(m_{(1)})]S^{-1}(1_{1})\\
   &=& 1_{1}'\c (1_{2}\c \a_M^{-1}(m_{(0)}))\o[\a_H^{-2}(1'_{2})\a_H^{-2}(m_{(1)})]S^{-1}(1_{1})\\
   &=& 1_{1}'\c (1_{2}\c \a_M^{-1}(m_{(0)}))\o\a_H^{-1}(1'_{2})[\a_H^{-2}(m_{(1)})S^{-1}(\a_H^{-1}(1_{1}))]\\
   &=& 1_{1}'\c (1_{2}\c \a_M^{-1}(m_{(0)}))\o1'_{2}[\a_H^{-2}(m_{(1)})S^{-1}(1_{1})].
   \end{eqnarray*}
Then we do a calculation as follows:
\begin{eqnarray*}
 &&(h_{2}\c m)_{(0)} \o (h_{2}\c m)_{(1)}\a_H^2(h_{1})\\
 &=&\a_H^{-1}(h_{221})\c m_{(0)}\o [(\a_H^{-2}(h_{222})\a_H^{-1}(m_{(1)}))S^{-1}(h_{21})]\a_H^2(h_{1})\\
 &=& h_{21}\c m_{(0)}\o [(\a_H^{-1}(h_{22})\a_H^{-1}(m_{(1)}))S^{-1}(h_{12})]\a_H(h_{11})\\
 &=&h_{21}\c m_{(0)}\o [(h_{22}m_{(1)})]S^{-1}(h_{12})h_{11}\\
 &=& h_{21}\c m_{(0)}\o [h_{22}m_{(1)}]S^{-1}(1_{1})\v(h_{1}1_{2})\\
&=&h_{12}\c m_{(0)}\o [\a_H(h_{2})m_{(1)}]S^{-1}(1_{1})\v(\a_H^{-1}(h_{11})1_{2})\\
&=&\a_H^{-1}(h_{12})1'_2\c m_{(0)}\o [\a_H(h_{2})m_{(1)}]S^{-1}(1_{1})\v(\a_H^{-2}(h_{11})\a_H^{-1}(1'_{1})1_{2})\\
&=&\a_H^{-1}(h_{12})1'_2\c m_{(0)}\o [\a_H(h_{2})m_{(1)}]S^{-1}(1_{1})\v(\a_H^{-1}(h_{11})[\a_H^{-1}(1'_{1})\a_H^{-1}(1_{2})])\\
&=&\a_H^{-1}(h_{12})1'_2\c m_{(0)}\o [\a_H(h_{2})m_{(1)}]S^{-1}(1_{1})\v(\a_H^{-1}(h_{11})[1'_{1}1_{2}])\\
 &=&h_{1}1_{2}\c m_{(0)}\o [h_{2}m_{(1)}]S^{-1}(1_{1})\\
 &=&h_{1}1'_{1}1_{2}\c m_{(0)}\o [\a_{H}(h_{2})[1'_{2}\a^{-1}_H(m_{(1)})]]S^{-1}(1_{1})\\
 &=&\a_H(h_{1})[1_{2}\c \a^{-1}_M(m_{(0)})]\o\a^{2}_{H}(h_{2})[[1_{3}\a^{-1}_H(m_{(1)})]S^{-1}(\a_{H}^{-1}(1_{1}))]\\
 &=&\a_H(h_{1})\c m_{(0)}\o \a^{2}_{H}(h_{2})\a_H(m_{(1)}).
\end{eqnarray*}
For $(3.2), (3.3)\Longrightarrow(3.1)$, we have
\begin{eqnarray*}
&&\a_H^{-1}(h_{21})\c m_{(0)}\o [\a_H^{-2}(h_{22})\a_H^{-1}(m_{(1)})]S^{-1}(h_{1})\\
&=& \a^{2}_M(\a_H^{-4}(h_{22})\c \a_M^{-2}(m))_{(0)}\o (\a_H^{-4}(h_{22})\c \a_M^{-2}(m))_{(1)}\a^{2}_H(\a_H^{-4}(h_{21}))S^{-1}(h_{1})\\
&=& (\a_H^{-1}(h_{2})\c m)_{(0)}\o (\a_H^{-2}(h_{2})\c \a_M^{-1}(m))_{(1)}\a_H^{-2}(h_{12}))S^{-1}(\a_H^{-2}(h_{11}))\\
&=&(\a_H^{-1}(h_{2})\c m)_{(0)}\o (\a_H^{-2}(h_{2})\c \a_M^{-1}(m))_{(1)}1_{1}\v(1_{2}\a_H^{-2}(h_{1}))\\
&=&(1_2'\a_H^{-2}(h_{2})\c m)_{(0)}\o (1_2'\a_H^{-3}(h_{2})\c \a_M^{-1}(m))_{(1)}1_{1}\v(1_{2}1_1'\a_H^{-3}(h_{1}))\\
&=&(1_3\a_H^{-2}(h_{2})\c m)_{(0)}\o (1_3\a_H^{-3}(h_{2})\c \a_M^{-1}(m))_{(1)}1_{1}\v(1_{2}\a_H^{-2}(h_{1}))\\
&=&(1_2\a_H^{-1}(h)\c m)_{(0)}\o (1_2\a_H^{-2}(h_{2})\c \a_M^{-1}(m))_{(1)}1_{1}\\
&=&(1_{2}\c(\a_H^{-1}(h)\c \a_H^{-1}(m))_{(0)}\o(1_{2}\c(\a_H^{-2}(h_{2})\c \a_M^{-2}(m)))_{(1)}1_{1}\\
&=&1_{1}\c\a_{M}^{-1}( h\c m)_{(0)}\o 1_{2}\a_{H}^{-1}(h\c m)_{(1)}\\
&=&(h\c m)_{(0)}\o (h\c m)_{(1)}.
\end{eqnarray*}

{\bf Definition 3.3.}  Let $(H,\a_H)$ be a weak Hom-Hopf algebra. Left-right weak-Hom type entwining structure is a triple $(A, C, \psi )$,
where $(A, \a_A)$ is a Hom-algebra and $(C, \a_C)$ is a  Hom-coalgebra with a
linear map  $\psi: A\o C\rightarrow A\o C$  such that $\psi\circ (\a_A\o \a_C)=(\a_A\o \a_C) \circ \psi$ satisfying the following conditions:
\begin{eqnarray}
{}_\psi(ab)\o \a_C(c^\psi)={}_\psi a {}_\varphi b\o \a_C(c)^{\psi\varphi},\\
\psi(c\o 1_A)=\varepsilon(c_1^{\psi}){}_\psi1_A\o c_2,\\
 \a_A({}_\psi a) \o \D(c^\psi)=\a_A(a)_{\varphi\psi})\o (c_{(1)}\ ^\psi \o c_{(2)}\ ^\varphi ),\\
\varepsilon(c^\psi){}_\psi a=\varepsilon(c^{\psi})a({}_\psi1_A).
\end{eqnarray}
Over a weak-Hom type entwining structure $(A, C, \psi)$, a left-right weak-Hom type entwined modules
$(M, \a_M)$ is both a right $C$-comodule and a left $A$-module such that
\begin{eqnarray}
\rho_{M}(a\c m )&=&\a_A(_\psi a) \c m_{(0)}\c \o \a_C(m_{(1)}^\psi),
\end{eqnarray}
for all $a\in A$ and $m\in M$.  ${}_{A}\cal{WM}(\psi)^{C}$ will denote  the category of left-right weak-Hom type entwined modules and morphisms between them.

{\bf Proposition 3.4.} Let $(H, \a_H)$ be a weak Hom-Hopf algebra. Define $\phi: H\o H\rightarrow H\o H$ given by $\phi(a\o c)=\a^{-1}_H(a_{21})\o (\a_H^{-2}(a_{22})\a_H^{-1}(c))S^{-1}(a_{1})$ for all $h,g\in H$, and so
${}_H{\cal WM}(\psi)^{H}$ is the category of  weak-Hom type entwined modules. In fact, for any $(M, \mu)\in {}_{H}{\cal WM}(\psi)^{H}$, one has compatible condition
\begin{eqnarray*}
\rho_M(a\c m)= \a^{-1}_H(a_{21}) \c m_{(0)}\o (\a_H^{-2}(a_{22})\a_H^{-1}(m_{(1)}))S^{-1}(a_{1}).
\end{eqnarray*}

{\it Proof. } We need to prove that  (3.4-3.7) hold. First, it is straightforward to check
(3.4) and (3.6). In what follows, we only verify (3.5) and (3.7). In fact, for all
$a,b,c\in H$, we have
\begin{eqnarray*}
 c_2\o \varepsilon(c_1^{\psi}){}_\psi1_A&=& c_2\o \varepsilon(\a_H^{-2}(1_{3})\a_H^{-1}(c_1))S^{-1}(1_{1}))\a_H(1_{2})\\
 &=& \widetilde{1}_2\a_H^{-1}(c_2)\o \varepsilon(1_{2}''[\widetilde{1}_1\a_H^{-2}(c_1))]1_2')\varepsilon(1_{1}'S^{-1}(1_{1}))1_{2}1_1''\\
 &=&\widetilde{1}_2\a_H^{-1}(c_2)\o \varepsilon([1_{2}''\widetilde{1}_1]\a_H^{-1}(c_1)1_2')\varepsilon(1_{1}'S^{-1}(1_{1}))1_{2}1_1''\\
 &=&1_3''\a_H^{-1}(c_2)\o \varepsilon([1_{2}''\a_H^{-1}(c_1)]1_2')\varepsilon(1_{1}'1_{2})S(1_{1})1_1''\\
 &=&1_3''\a_H^{-2}(c_2)1_3'\o \varepsilon([1_{2}''\a_H^{-1}(c_1)]1_2')\varepsilon(1_{1}'1_{2})S(1_{1})1_1''\\
 &=&1_2''\a_H^{-1}(c)1_2'\o \varepsilon(1_{1}'1_{2})S(1_{1})1_1''\\
 &=&1_2''\a_H^{-1}(c)1_2\o S(1_{1})1_1''\\
 &=&1_2'\a_H^{-1}(c)S^{-1}(1_{1})\o 1_{2}1_1'\\
  &=&1_3\a_H^{-1}(c)S^{-1}(1_{1})\o 1_{2}\\
  &=& c^{\psi}\o {}_\psi 1_A.
\end{eqnarray*}
As for (3.7), we compute:
\begin{eqnarray*}
\varepsilon(c^\psi){}_\psi a&=& \varepsilon((\a_H^{-2}(a_{22})\a_H^{-1}(c))S^{-1}(a_{1})) \a_H(a_{21})\\
&=& \varepsilon((\a_H^{-2}(a_{22})\a_H^{-1}(c))1_2) \varepsilon(1_1S^{-1}(a_{1})) \a_H(a_{21})\\
&=&\varepsilon(\a_H^{-1}(a_{22})(\a_H^{-1}(c)1_2)) \varepsilon(1_1S^{-1}(a_{1})) \a_H(a_{21})\\
&=&\varepsilon(\a_H^{-1}(a_{22})1_1') \varepsilon(1_2'(\a_H^{-1}(c)1_2))\varepsilon(1_1 S^{-1}(a_{1})) \a_H(a_{21})\\
&=&\varepsilon(\a_H^{-1}(a_{22})1_1') \varepsilon(1_2'(\a_H^{-1}(c)S^{-1}(1_1)))\varepsilon(a_{1}1_2) \a_H(a_{21})\\
&=&\varepsilon(a_{2}1_1') \varepsilon(1_2'(\a_H^{-1}(c)S^{-1}(1_1)))\varepsilon(\a_H^{-1}(a_{11})1_2) \a_H(a_{12})\\
&=&\varepsilon(a_{2}1_1') \varepsilon(1_2'(\a_H^{-1}(c)S^{-1}(1_1)))\varepsilon(a_{1}(1_1''1_2)) a_{21}1_2''\\
&=&\varepsilon(\a^{-1}(a_{2})1_2''1_1') \varepsilon(1_2'(\a_H^{-1}(c)S^{-1}(1_1)))a_{1}1_1''1_2\\
&=&\varepsilon(a_{2}[1_2''1_1']) \varepsilon(1_2'(\a_H^{-1}(c)S^{-1}(1_1)))a_{1}1_1''1_2\\
&=& \varepsilon(1_2'(\a_H^{-1}(c)S^{-1}(1_1)))\a_H(a_{1})1_1'1_2\\
&=&\varepsilon(1_2'(\a_H^{-1}(c)S^{-1}(1_1)))a_{1}[1_1'1_2]\\
&=&\varepsilon(1_3(\a_H^{-1}(c)S^{-1}(1_1)))a_{1}1_2\\
&=&\varepsilon(c^{\psi})a({}_\psi1_A).
\end{eqnarray*}

{\bf Proposition 3.5.} Let $(H,\a_H)$ be a weak Hom-Hopf algebra, for any $(M,\a_M),  (N,\a_N)  \in
{}_H{\cal WYD}^H$, and define the linear map
\begin{eqnarray*}
B_{M,N}: M\o N\rightarrow N\o M,~~~~~~~~B_{M,N}(m\o n)=n_{(0)}\o \a^{-1}_H(n_{(1)})\c m.
\end{eqnarray*}
Then, we have $(\a_N \o \a_M)\circ B_{M,N} = B_{M,N}\circ (\a_M\o \a_N)$ and, if $(P, \a_P)\in {}_H{\cal WYD}^H$, the maps $B_{-,-}$ satisfy the Hom-Yang-Baxter
 equation:
 \begin{eqnarray*}
&&(\a_P\o B_{M,N})\circ (B_{M,P})\o (\a_M\o B_{N,P})\\
&=&(B_{N,P}\o \a_M)\circ ( \a_N \o B_{N,P})\circ ((B_{M,P}\o \a_M).
 \end{eqnarray*}
{\it Proof.} The proof is similar to Proposition 3.4 in \cite{MP14}.

{\bf lemma 3.6.} Let$(H,\a_H)$ be a weak Hom-Hopf algebra, then $H_s$ is the unit object in ${}_H{\cal WYD}^H$ with the action:   for any $h \in H$, $x\in H_s$,
$$
h\cdot x = \widehat{\varepsilon_s}(h)(hx),~~ \rho(x)=x_1\o x_2,
$$
and $\a_{H_s} = \a$.

{\it Proof.} The proof is similar to \cite{N02}.

{\bf lemma 3.7.} Let $(H,\a_H)$ be a weak Hom-Hopf algebra,  the left and right unit constraints $l_M: H_s\o_t M\rightarrow M$ and $r_M: M\o_t H_s\rightarrow M$ and
their inverses are given by the formulas
\begin{eqnarray*}
&&l_M(x \o m)=S(x)\c\a^{-2}_M(m),~~~~~~~~~~~l_M^{-1}(m)=1_H\o  \a_M(m),\\
&&r_M(m\o x)=x\c \a^{-2}_M(m),~~~~~~r_M^{-1}(m)=\varepsilon(1_3)\varepsilon_s(1_2)\c \a_M(m)\o 1_1.
\end{eqnarray*}
{\it Proof.} It is easy to see that $l_M$ is natural isomorphisms in ${}_H{\cal WYD}^H$.
 We only check that
$$\aligned
l_{M}^{-1}l_{M}(x \o_t m)=~& l_{M}^{-1}(S(x)\cdot \a_M^{-2}(m)) = 1_H \o_t S(x)\cdot \a_M^{-1}(m)\\
=~&\varepsilon_t (1_{1})\o (1_{2}S(x))\cdot m\\
=~&\varepsilon_t ({\widehat{\varepsilon_s} (S(x))}_{1}) \o {\widehat{\varepsilon_s} (x)}_{2}\cdot m\\
=~&\varepsilon_t (1_{1}\widehat{\varepsilon_s} (S(x))) \o 1_{2}\cdot m \\
=~&\varepsilon_t (1_{1}x) \o 1_{2}\cdot m = x \o_t m,
\endaligned$$
and
$$\aligned
l_{M}l_{M}^{-1}(m) = l_{M}( 1_H\o_t \a_M(m)) = 1_H\cdot \a^{-2}_M(\a_M(m)) = m,
\endaligned$$
which implies $l^{-1}_M$ is the inverse of $l_M$.

Similarly, we can check that $r_M$ is a natural isomorphism with the inverse $r^{-1}_M$ in ${}_H{\cal WYD}^H$.

{\bf Theorem 3.8.}  Let $(H,\a_H)$ be a weak Hom-Hopf algebra. Then $({}_H{\cal WYD}^H,\o_t, H_s)$ is a monoidal category.

{\it Proof.} Firstly, for any $(M, \a_M),(N,\a_N), (P,\a_P) \in {}_H{\cal WYD}^H$, define an associativity constraint by
$$
a_{M,N,P}((m \o_t n) \o_t p) = \a^{-1}_M(m) \o_t(n \o_t \a_P(p)),~m \in M,~n \in N,~p \in P.
$$
Obviously that $a$ is natural and satisfies $a_{M,N,P}\circ(\a_M \o (\a_N \o \a_P)) = ((\a_M \o \a_N)\o \a_P)\circ a_{M,N,P}$. For any $h \in H$, since
$$\aligned
&~~~~a_{M,N,P}(h\cdot ((m \o_t n) \o_t p))\\
&= \a_M^{-1}(h_{11}\cdot m )\o_t(h_{12}\cdot n \o_t \a_P(h_{2}\cdot p))\\
&=h_1\cdot \a_M^{-1}(m) \o_t( h_{21}\cdot n \o_t h_{22}\cdot \a_P(p))\\
&=h\cdot(\a_M^{-1}(m) \o_t(n \o_t \a_P(p) ))\\
&= h\cdot (a_{M,N,P}((m \o_t n) \o_t p) ),
\endaligned$$
$a_{M,N,P}$ is $H$-linear.

Next we will check that $a_{M,N,P}$ is $H$-colinear.

$$\aligned
&(a_{M,N,P}\o id_H)\circ \rho_{(M \hat{\o}N)\hat{\o}P}((m\o n)\o p)\\
=~& (a_{M,N,P}\o id_H)((m\o n)_{(0)} \o p_{(0)}\o \a^{-2}_H(p_{(1)}(m\o n)_{(1)}))\\
=~&a_{M,N,P}((m_{(0)}\o n_{(0)})\o p_{(0)})\o  \a^{-2}_H(p_{(1)}\a^{-2}_H(n_{(1)}m_{(1)}))\\
=~&\a_M^{-1}(m_{(0)})\o (n_{(0)}\o \a_P(p_{(0)}))\o \a^{-2}_H(p_{(1)})\a^{-4}_H(n_{(1)}m_{(1)}),
\endaligned$$
$$\aligned
&\rho_{(M \hat{\o}N)\hat{\o}P}\circ a_{M,N,P}((m\o n)\o p)\\
=~& \rho_{(M \hat{\o}N)\hat{\o}P}(\a_M^{-1}(m)\o (n\o \a_P(p)))\\
=~& \a_M^{-1}(m)_{(0)}\o (n\o \a_P(p))_{(0)}\o \a_H^{-2}((n\o \a_P(p))_{(1)}\a_M^{-1}(m)_{(1)})\\
=~& \a_M^{-1}(m)_{(0)}\o (n_{(0)}\o \a_P(p)_{(0)})\o \a_H^{-2}(\a^{-2}_H(\a_H(p)_{(1)} n_{(1)})\a_M^{-1}(m)_{(1)})\\
=~& \a_M^{-1}(m)_{(0)}\o (n_{(0)}\o \a_P(p)_{(0)})\o [\a_H^{-3}(p)_{(1)} \a_H^{-4}(n_{(1)})]\a_M^{-3}(m)_{(1)}\\
=~&\a_M^{-1}(m_{(0)})\o (n_{(0)}\o \a_P(p_{(0)}))\o [\a_H^{-2}(p)_{(1)} [\a_H^{-4}(n_{(1)})\a_M^{-4}(m)_{(1)}]\\
=~&\a_M^{-1}(m_{(0)})\o (n_{(0)}\o \a_P(p_{(0)}))\o \a^{-2}_H(p_{(1)})\a^{-4}_H(n_{(1)}m_{(1)}).
\endaligned$$
And $a_{M,N,P}$ is a bijection because of $\a_M$, $\a_P$ are all bijective maps. Thus $a$ is a natural isomorphism in ${}_H{\cal WYD}^H$.

Secondly, it is also a direct check to prove that $a$ satisfies the Pentagon Axiom.

At last, we will check the Triangle Axiom. In fact, for any $x \in H_s$, we have
$$\aligned
&(r_M \o id_N)((m \o_t x) \o_t n)\\
=~& x\cdot \a^{-2}_M(m) \o_t n \\
=~&1_1\c (x\cdot \a^{-2}_M(m)) \o_t 1_2\c n \\
=~&1_1x\c \a^{-1}_M(m)\o_t 1_2\c n \\
=~&1_1\c \a^{-1}_M(m)\o_t 1_2S(x)\c n \\
=~& 1_1\c \a^{-1}_M(m)\o_t 1_2\c (S(x)\c \a^{-1}_N(n) )\\
=~& \a^{-1}_M(m) \o S(x)\c \a^{-1}_N(n)\\
=~& \a^{-1}_M(m) \o_t S(x)\cdot\a_N^{-2}(\a_N(n))\\
=~& (id_M \o l_N)a_{M,H_t ,N}((m \o_t x) \o_t n).
\endaligned$$

Let $(H,\a_H)$ be a weak Hom-Hopf algebra with a bijective antipode $S$. Consider the full subcategory ${}_H{\cal WYD}^H_{f.d.}$ of ${}_H{\cal WYD}^H$ whose objects are  finite-dimensional. Using the antipode $S$ of $H$, we can provide ${\cal WYD}(H)_{f.d.}$ with a duality.

For any $(M, \a_M) \in {}_H{\cal WYD}^H_{f.d.}$, we set ${}^{\ast}M=Hom(M,k)$,  with the action and the coaction of $H$ on
$M^\ast$ given by
$$(h\c f)(m)=f(S(\a^{-1}_H(h))\c \a_M^{-2}(m))~~~ \mbox{and}~~~  f_{(0)}(m)\o
f_{(1)}=f(\a_M^{-2}(m_{(0)}))\o S^{-1}(\a^{-1}_H(m_{(1)})).$$

Similarly,  for any  $(M, \a_M) \in {\cal WYD}(H)_{f.d.}$,  we set ${}^{\ast}M=Hom(M,k)$,   with the action and the coaction of $H$ on
$M^\ast$ given by
$$(h\c f)(m)=f(S^{-1}(\a^{-1}_H(h))\c \a_M^{-2}(m))~~~ \mbox{and}~~~  f_{(0)}(m)\o
f_{(1)}=f(\a_M^{-2}(m_{(0)}))\o S(\a^{-1}_H(m_{(1)})).$$

{\bf Theorem 3.9.}  The category ${}_H{\cal WYD}^H_{f.d.}$ is a rigid category.

{\it Proof.}  Define the maps
$$
coev_M: ~H_t \rightarrow M \o_t M^\ast,~~~~x\mapsto \sum x\cdot(e_i \o_t \a_{M^\ast}(e^i)),
$$
where $e_i$ and $e^i$ are bases of $M$ and $M^\ast$, respectively, dual to each other, and
$$
ev_M:~M^\ast \o_t M \rightarrow H_t,~~~~f \o_t m \mapsto f(1_{1}\cdot m)1_{2}.
$$

Firstly, we will prove that $M^{*}$ is indeed an object in ${}_H{\cal WYD}^H$,  and $\a_{{}^\ast M}$ is given by
$$
\a_{{}^\ast M}(f)(m) = f(\a^{-1}_M(m)),~~~~f \in {}^\ast M,~m \in M.
$$
We have
\begin{eqnarray*}
(h\c f)_{(0)}(m)\o (h\c f)_{(1)}&=&(h\c f)(\a^{-2}_M(m_{(0)}))\o
S^{-1}(\a^{-1}_H(m_{(1)}))\\
 &=&f(S(\a^{-1}_H(h))\c \a^{-4}(m_{(0)})))\o S^{-1}(\a^{-1}_H(m_{(1)})),
 \end{eqnarray*}
 \begin{eqnarray*}
 &&(\a^{-1}_H(h_{21})\c f_{(0)})(m)\o (\a^{-2}_H(h_{22})\a^{-1}_H(f_{(1)})) S^{-1}(h_{1})\\
 &=&f_{(0)}(S(\a^{-2}_H(h_{21}))\c \a^{-2}_M(m))\o (\a^{-2}_H(h_{22})\a^{-1}_H(f_{(1)}))S^{-1}(h_{1})\\
 &=&f(S(\a^{-4}_H(h_{21}))\c \a^{-4}_M(m))_{(0)})\o (\a^{-2}_H(h_{22})\a^{-2}_H(S(\a^{-1}_H(h_{21}))\c \a^{-2}_H(m_{(1)})))S^{-1}(h_{1})\\
 &=&f(S(\a^{-5}_H(h_{2112}))\c \a^{-4}_M(m)_{(0)})\o (\a^{-2}_H(h_{22})\\
 &&[S^{-1}(\a^{-4}_H(h_{212}))(\a^{-5}_H(S^{-1}(m_{(1)}))\a^{-6}_H(h_{2111}))])S^{-1}(h_{1})\\
 &=&f(S(\a^{-4}_H(h_{212}))\c \a^{-4}_M(m)_{(0)}))\o ([\a^{-4}_H(h_{222})S^{-1}(\a^{-4}_H(h_{221}))]\\
 &&(\a^{-4}_H(S^{-1}(m_{(1)}))\a^{-4}_H(h_{211}))])S^{-1}(h_{1})\\
 &=&f(S(\a^{-4}_H(h_{212}))\c \a^{-2}_M(m)_{(0)}))\o (1_1\varepsilon(1_2h_{22})\\
 &&(\a^{-4}_H(S^{-1}(m_{(1)}))\a^{-4}_H(h_{211}))]) S^{-1}(h_{1})\\
  &=&f(S(\a^{-4}_H(h_{221}))\c \a^{-4}_M(m)_{(0)}))\o (1_1\varepsilon(1_2h_{222})(\a^{-4}_H(S^{-1}(m_{(1)}))\a^{-3}_H(h_{21}))])S^{-1}(h_{1})\\
   &=&f(S(\a^{-4}_H(h_{221}))\c \a^{-4}_M(m)_{(0)}))\o (1_1\varepsilon(1_2h_{222})(\a^{-4}_H(S^{-1}(m_{(1)}))\a^{-3}_H(h_{21}))])S^{-1}(h_{1})\\
 &=&f(S(\a^{-3}_H(h_{21}))\c \a^{-4}_M(m)_{(0)}))\o (1_1\varepsilon(1_2h_{22})(\a^{-3}_H(S^{-1}(m_{(1)}))[\a^{-2}_H(h_{12}))S^{-1}(\a^{-2}_H(h_{11}))]\\
 &=&f(S(\a^{-3}_H(h_{21}))\c \a^{-4}_M(m)_{(0)}))\o (1_1'\varepsilon(1'_2h_{22})(\a^{-3}_H(S^{-1}(m_{(1)}))1_1\varepsilon(1_2h_1)\\
  &=&f(S(\a^{-4}_H(1_2''h_{12}))\c \a^{-4}_M(m)_{(0)}))\o (1_1'\varepsilon(1'_2h_{2})(\a^{-2}_H(S^{-1}(m_{(1)}))1_1\varepsilon(1_21''_1h_{11})\\
   &=&f(S(\a^{-3}_H(1_2h_{1}))\c \a^{-4}_M(m)_{(0)}))\o (1_1'\varepsilon(1_2'h_{2})(\a^{-3}_H(S^{-1}(m_{(1)}))1_1\\
   &=&f(S(\a^{-4}_H(1''_{1}h_{1}))S(1_2)\c \a^{-4}_M(m)_{(0)}))\o (S^{-1}(1'_{2})\varepsilon(1_1'1_2''h_{2})(\a^{-3}_H(S^{-1}(m_{(1)}))1_1\\
  &=&f(S(\a^{-3}_H(1_21'_{1}h))\c \a^{-4}_M(m)_{(0)}))\o (S^{-1}(1'_{2})(\a^{-3}_H(S^{-1}(m_{(1)}))1_1\\
 &=&f(S(\a^{-2}_H(h))S(1_2)\c \a^{-4}_M(m)_{(0)}))\o (S^{-1}(1_3\a^{-3}_H(S^{-1}(m_{(1)}))S^{-1}(1_1))\\
 &=&f(S(\a^{-1}_H(h))(1_2\c \a^{-5}_M(m)_{(0)}))\o (S^{-1}(1_3\a^{-3}_H(S^{-1}(m_{(1)}))S^{-1}(1_1))\\
 &=&f(S(\a^{-1}_H(h))\c \a^{-4}_M(m)_{(0)})\o S^{-1}(\a^{-1}_H(m_{(1)})).
  \end{eqnarray*}
 Which means that
 $$(h\c f)_{(0)}\o (h\c f)_{(1)}=(\a^{-1}_H(h_{21})\c f_{(0)})\o
 (\a_H^{-2}(h_{22})\a_{H}^{-1}(f_{(1)}))S^{-1}(h_{1}).$$

We have known that $H_s\in {}_H{\cal WYD}^H$, with left $H$-module structure $h\c z = \widehat{\varepsilon}_s(hz)$ and
right $H$-comodule structure $\rho(x)=x_1\o x_2$, for all $x\in H_s$. Next,  we will check $ev_M$ and $coev_M$ are morphisms in ${}_H{\cal WYD}^H$. For any $h \in H$, $m \in M$, $f \in M^\ast$, we compute

 \begin{eqnarray*}
ev_M(h\cdot(f \o_t m))&=& (h_{1}\cdot f)(1_{1}\cdot(h_{2}\cdot m))1_{2} \\
&=&f(S(1'_1\a_H^{-2}(h_1))\cdot((1_1\cdot(1'_2\a_H^{-3}(h_2)))\cdot\a_{M}^{-1}(m)))1_2 \\
&=&f((S(1_1\a_H^{-2}(h_1))(1_2\a_H^{-2}(h_2)))\cdot\cdot\a_{M}^{-1}(m))1_3 \\
&=&f(\varepsilon_s (1_1 \a_H^{-2}(h))\cdot m)1_2\\
&=&f(\varepsilon_s (\a_H^{-2}(h_1))\cdot v)\varepsilon_t(\a_H^{-2}(h_1)) \\
&=&f(1_1\cdot m)\varepsilon_t(h1_2) = h\cdot (ev_M(f \o_t m)),
 \end{eqnarray*}
and
on one hand, \begin{eqnarray*}
\rho \circ ev_M(f \o m)&=& \rho(1_2)f(1_1\c m)\\
&=&f(1_1\c m)S(1_1')\o 1_2'S^{-1}(1_2),
\end{eqnarray*}
on the other hand
\begin{eqnarray*}
(ev_M\o id)\circ\rho(f\o m)&=& (ev_M\o id)(f_{(0)}\o m_{(0)})\o \a^{-2}_H(m_{(1)}f_{(1)})\\
&=& f_{(0)}(1_1\c m_{(0)})1_2\o \a^{-2}_H(m_{(1)}f_{(1)})\\
&=& f(1_2\c \a_{M}^{-2}(m_{(0)(0)}))(1_2')\o \a^{-2}_H(m_{(1)})S^{-1}(1_1'1_3\a^{-3}_H(m_{(0)(1)})S^{-1}(1_1))\\
&=& f(1_2\c \a_{M}^{-1}(m)_{(0)})(1_2')\o \a^{-3}_H(m_{(1)2})S^{-1}(1_1'1_3\a^{-3}_H(m_{(1)1})S^{-1}(1_1))\\
&=&f(1_2\c \a_{M}^{-1}(m)_{(0)})(1_2')\o S^{-1}(\varepsilon_t(m_{(1)}))S^{-1}(1_1'))\\
&=&f(1_1\c m)1_2'\o S^{-1}(1_2)S^{-1}(1_1')\\
&=&f(1_1\c m)S(1_1')\o 1_2'S^{-1}(1_2).
\end{eqnarray*}
Thus $ev_M$ is $H$-linear and $H$-colinear. And it is easy to get that $ev_M\circ(\a_{M^\ast} \o \a_M) = \a_H\circ ev_M$, hence $ev_M \in {}_H{\cal WYD}^H$.

Next we have
 \begin{eqnarray*}
coev_M(h\cdot x)(m) &=&\sum \varepsilon_t (hx)\cdot(e_i \o \a_{M^\ast}(e^i))(m)\\
&=&\varepsilon_t (hx)\cdot \a^{-2}_M(m)),
 \end{eqnarray*}
and
 \begin{eqnarray*}
h\cdot coev_M(x)(m) &=& (\a_H^{-1}(h_1)x_1)\cdot\a_M(e_i) \o_t ((\a_H^{-1}(h_2)x_2)\cdot\a^2_{M^\ast}(e^i))(m)\\
&=&(\a_H^{-1}(h_1)x_1)\cdot(S(\a_H^{-2}(h_2)\a_H^{-1}(x_2))\cdot \a^{-3}_M(m))\\
&=&\varepsilon_t(\a_H^{-2}(h)\a_H^{-1}(x))\cdot \a^{-2}_M(m)=\varepsilon_t (hx)\cdot \a^{-2}_M(m),
 \end{eqnarray*}
hence $coev_V$ is $H$-linear,  it is not hard to check that $coev_V$ is $H$-colinear and is left to the reader.   Obviously that $coev_M\circ\a_H = (\a_M \o \a_{M^\ast})\circ coev_M$, thus $coev_M \in {}_H{\cal WYD}^H$.

Secondly, we consider
 \begin{eqnarray*}
&&(r_M\circ(id_M \o ev_M)\circ a_{M,M^\ast,M} \circ (coev_M \o id_M)\circ l_M^{-1})(m)\\
&=&(r_M\circ(id_M \o ev_M)\circ a_{M,M^\ast,M})((1_1\cdot e_i \o_t 1_{2}\cdot \a_{M^\ast}(e^i)) \o_t \a_M(m))\\
&=&r_M(e_i \o_t 1_2\cdot\a_{M^\ast}^2(e^i) (1_1\cdot\a^2_M(m)) )\\
&=&\widehat{\varepsilon_s}(1_2)\cdot(1_1\cdot \a^2_M(m))= m,
 \end{eqnarray*}
and
 \begin{eqnarray*}
&&(l_{M^\ast}\circ (ev_M \o id_{M^\ast})\circ a^{-1}_{M^\ast,M,M^\ast}\circ (id_{M^\ast} \o coev_M)\circ r_{M^\ast}^{-1})(f)(m)\\
&=&(l_{M^\ast}\circ (ev_M \o id_{M^\ast})\circ a^{-1}_{M^\ast,M,M^\ast})(\a_{M^\ast}(f) \o_t (\a_{M }(e_i) \o_t \a^2_{M^\ast}(e^i)))(m) \\
&=&l_{V^\ast}(\a^2_{M^\ast}(f)(1_1\cdot \a_{M }(e_i)) 1_2 \o_t\a_{M^\ast}(e^i))(m) \\
&=&f(1_1\cdot(S(1_2)\cdot \a^{-2}_M(m))) = f(m).
 \end{eqnarray*}
Thus ${}_H{\cal WYD}^H_{f.d.}$ is a left rigid category.

Similarly, we define the following maps

$$
\widetilde{coev}_M: ~H_t \rightarrow {}^\ast M \o_t M,~~~~x\mapsto x\cdot(\sum \a_{{}^\ast M}(e^i) \o_t e_i),
$$
and
$$
\widetilde{ev}_M:~M \o_t {}^\ast M \rightarrow H_t ,~~~~ m \o_t f\mapsto f(S^{-1} (1_{1}) \cdot m)1_{2}.
$$
We can show that $({}^\ast M,\widetilde{ev}_M, \widetilde{coev}_M)$ is a right dual of $M$.
Thus ${}_H{\cal WYD}^H_{f.d.}$ is a right rigid category.

\section*{4. A braided monoidal category ${}_H{\cal WYD}^H$ I}
\def\theequation{4. \arabic{equation}}
\setcounter{equation} {0} \hskip\parindent

{\bf Proposition 4.1.} Let $(H,\a_H)$ be a weak Hom-Hopf algebra. For any $(M,\a_M),  (N,\a_N)  \in
{}_H{\cal WYD}^H$,   then  $M \o_{t} N=1_{1}M \o 1_{2}N \in
{}_H{\cal WYD}^H$ with structures:
\begin{eqnarray*}
&&h\c (m \o_t n)=h_{1}\c m \o_t h_{2}\c n,\\
&&m\o_t n \mapsto(m\o_t n)_{(0)}\o_t (m\o_t n)_{(1)}=(m_{(0)}\o_t n_{(0)})\o
\a^{-2}_H(n_{(1)}m_{(1)}).
\end{eqnarray*}
{\it Proof.} Obviously $M\o_{t}N$ is a left $H$-module and a right
$H$-comodule. We check now the compatibility condition. We compute:
\begin{eqnarray*}
&&(h\c (m\o n))_{(0)}\o (h\c (m\o n))_{(1)}\\
&=&(h_{1}\c m \o h_{2}\c n)_{(0)}\o (h_{1}\c m
\o h_{2}\c n)_{(1)}\\
&=&((h_{1}\c m)_{(0)} \o (h_{2}\c n)_{(0)}\o
\a^{-2}_H ((h_{2}\c n)_{(1)} (h_{1}\c m)_{(1)})\\
&=& (\a_H^{-1}(h_{121})\c m_{(0)}\o \a_H^{-1}(h_{221})\c n_{(0)})\o \a^{-2}_{H}(((\a^{-2}_H(h_{222})\a^{-1}_H(n_{1}))S^{-1}(h_{21}))\\
&&((\a^{-2}_H(h_{122})\a^{-1}_H(m_{(1)}))S^{-1}(h_{11})))\\
&=& (h_{12}\c m_{(0)}\o \a_H^{-1}(h_{212})\c n_{(0)})\o \a^{-2}_{H}((h_{22}n_{1})\\
&&((S^{-1}(\a^{-3}_H(h_{2112}))\a^{-3}_H(h_{2111}))(\a^{-1}_H(m_{(1)})S^{-1}(\a^{-1}_H(h_{11}))))\\
&=& (h_{12}\c m_{(0)}\o \a_H^{-1}(h_{212})\c n_{(0)})\o \a^{-2}_{H}((h_{22}n_{1})\\
&&(1_2\varepsilon(\a^{-3}_H(h_{211})1_1)(\a^{-1}_H(m_{(1)})S^{-1}(\a^{-1}_H(h_{11}))))\\
&=& (h_{12}\c m_{(0)}\o \a_H^{-1}(h_{221})\c n_{(0)})\o \a^{-2}_{H}((h_{222}n_{1})\\
&&(1_2\varepsilon(\a^{-2}_H(h_{21})1_1)(\a^{-1}_H(m_{(1)})S^{-1}(\a^{-1}_H(h_{11}))))
\end{eqnarray*}
\begin{eqnarray*}
&=& (\a^{-1}_H(h_{211})\c m_{(0)}\o \a_H^{-1}(h_{221})\c n_{(0)})\o \a^{-2}_{H}((h_{222}n_{1})\\
&&(1_2\varepsilon(\a^{-2}_H(h_{212})1_1)(\a^{-1}_H(m_{(1)})S^{-1}(h_{1})))\\
&=& (\a^{-2}_H(h_{211})1_1'\c m_{(0)}\o \a_H^{-1}(h_{221})\c n_{(0)})\o \a^{-2}_{H}((h_{222}n_{1})\\
&&(1_2\varepsilon(\a^{-2}_H(h_{212})1_2'1_1)(\a^{-1}_H(m_{(1)})S^{-1}(h_{1})))\\
&=& (\a^{-2}_H(h_{211})1_1'\c m_{(0)}\o \a_H^{-1}(h_{221})\c n_{(0)})\o \a^{-2}_{H}((h_{222}n_{1})\\
&&(1_2\varepsilon(\a^{-2}_H(h_{212})1_2'1_1)(\a^{-1}_H(m_{(1)})S^{-1}(h_{1})))\\
&=& (\a^{-2}_H(h_{211})1_2'1_1\c m_{(0)}\o \a_H^{-1}(h_{212})\c n_{(0)})\o \a^{-2}_{H}((\a_H(h_{22})n_{1})\\
&&(1_2(\a^{-1}_H(m_{(1)})S^{-1}(\a^{-1}(h_{1})1_{1}')))\\
&=&\a^{-1}_H(h_{21})\c (m\o n)_{(0)}\o\a^{-2}_{H}(h_{22})\a^{-1}_H(m\o
n)_{(1)}S^{-1}(h_{1}).
\end{eqnarray*}
Hence $ M \o_{t} N \in {}_H{\cal WYD}^H$.

{\bf Proposition 4.2.} Let $(M,\a_M), (N,\a_N)\in {}_H{\cal WYD}^H$.  Define the map
$$ c_{M,N}: M \o_{t} N \rightarrow M \o_{t} M,\,\,\, c_{M,N}(m\o
n)=\a^{-1}_{N}(n_{(0)})\o \a_M^{-1}(\a^{-1}_H(n_{(1)})\c m).$$
 Then $ c_{M,N}$ is $H$-linear
$H$-colinear and satisfies the conditions (for $(P, \a_P) \in {}_H{\cal
WYD}^H$)
\begin{equation}
a^{-1}_{P, M,N}\circ c_{M\o N,P}\circ a^{-1}_{M,N,P}=(c_{M,P}\o id_{N})\circ a^{-1}_{M,P,N}\circ (id_{M}\o c_{N,P}),
\end{equation}
\begin{equation}
a_{N,P,M}\circ c_{M,N\o P}\circ a_{M,N,P}=(id_{N} \o c_{M,P} )\circ a_{N,M,P}\circ (c_{M,N}\o id_{P}).
\end{equation}

{\it Proof.}  First, we prove that $c_{M,N}$ is $H$-linear, we compute:
\begin{eqnarray*}
&& c_{M,N}(h\c (m\o n))\\
&=&c_{M,N} (h_{1}\c m\o h_{2}\c n)\\
&=&\a^{-1}_{N}((h_2\c n)_{(0)})\o \a_M^{-1}(\a^{-1}_H((h_2\c n)_{(1)}))\c (h_{1}\c m))\\
&=& \a^{-1}_{N}(\a^{-1}_H(h_{221})\c n_{(0)})\o \a_M^{-1}(\a^{-1}_H([\a^{-2}_H(h_{222})\a^{-1}_H(n_{(1)})]S^{-1}(h_{21}))\c (h_{1}\c m))\\
&=& \a^{-1}_{N}(h_{21}\c n_{(0)})\o \a_M^{-1}(\a^{-1}_H([\a^{-1}_H(h_{22})\a^{-1}_H(n_{(1)})]S^{-1}(h_{12}))\c (\a_{H}^{-1}(h_{11})\c m))\\
&=&\a^{-1}_{N}(h_{21}\c n_{(0)})\o \a_M^{-1}([\a^{-2}_H(h_{22})\a^{-2}_H(n_{(1)})]S^{-1}(\a_{H}^{-1}(h_{12}))\c (\a_{H}^{-1}(h_{11})\c m))\\
&=&\a^{-1}_{N}(h_{21}\c n_{(0)})\o \a_M^{-1}([\a^{-2}_H(h_{22})\a^{-2}_H(n_{(1)})][S^{-1}(\a_{H}^{-2}(h_{12}))\a_{H}^{-2}(h_{11})]\c \a_M(m))\\
&=&\a^{-1}_{N}(h_{21}\c n_{(0)})\o \a_M^{-1}([\a^{-2}_H(h_{22})\a^{-2}_H(n_{(1)})][S^{-1}(1_1)\varepsilon(\a_{H}^{-2}(h_{1}1_2))]\c \a_M(m))\\
&=&\a^{-1}_{N}(\a^{-1}_H(h_{12})1_3\c n_{(0)})\o \a_M^{-1}([\a^{-2}_H(h_{2})1_4\a^{-2}_H(n_{(1)})][S^{-1}(1_1)\varepsilon(\a_{H}^{-2}(h_{11}1_2))]\c \a_M(m))\\
&=&\a^{-1}_{N}(h_{1}1_2\c n_{(0)})\o \a_M^{-1}([\a^{-2}_H(h_{2})1_3\a^{-2}_H(n_{(1)})]S^{-1}(1_1)\c \a_M(m))\\
&=&h_{1}\c (1_2\c \a^{-2}_{N}(n_{(0)}))\o \a^{-1}_H(h_{2})[1_3\a^{-4}_H(n_{(1)})S^{-1}(1_1)]\c m\\
&=& h_1 \c \a^{-1}_{N}(n_{(0)})\o \a_M^{-1}(h_2)\a^{-2}_H(n_{(1)})\c m\\
&=& h\c c_{M,N}(m\o n).
\end{eqnarray*}
Next we prove that $c_{M,N}$ is $H$-colinear.
\begin{eqnarray*}
&& \rho_{N\o M}c_{M,N}(m\o n)\\
&=& \rho_{N\o M}(\a^{-1}_{N}(n_{(0)})\o \a_M^{-1}(\a^{-1}_H(n_{(1)})\c m))\\
&=&\a^{-1}_{N}(n_{(0)(0)})\o (\a^{-2}_H(n_{(1)})\c \a_M^{-1}(m))_{(0)}\o \a_{H}^{-2}((\a^{-2}_H(n_{(1)})\c \a_M^{-1}(m))_{(1)}\a^{-1}_{N}(n_{(0)(1)}))\\
&=& \a^{-1}_{N}(n_{(0)(0)})\o \a^{-3}_H(n_{(1)21})\c \a_M^{-1}(m_{(0)})\o \a_{H}^{-2}((\a^{-4}_H(n_{(1)22}) \a_H^{-2}(m_{(1)})\\
&&S^{-1}(n_{(1)1}))\a^{-1}_{N}(n_{(0)(1)}))\\
&=&\a^{-1}_{N}(n_{(0)(0)})\o \a^{-2}_H(n_{(1)1})\c \a_M^{-1}(m_{(0)})\o \a_{H}^{-2}((\a^{-2}_H(n_{(1)2}) \a_H^{-1}(m_{(1)}))\\
&&[S^{-1}(\a^{-3}_H(n_{(0)(1)2}))\a^{-3}_{N}(n_{(0)(1)(1)})])\\
&=&n_{(0)}\o \a^{-3}_H(n_{(1)1})1_2\c \a_M^{-1}(m_{(0)})\o \a_{H}^{-2}((\a^{-3}_H(n_{(1)2})1_3 \a_H^{-1}(m_{(1)}))S^{-1}(1_1))\\
&=& \a^{-1}_{N}(n_{(0)(0)})\o \a^{-3}_H(n_{(0)(1)})1_2\c \a_M^{-1}(m_{(0)})\o \a_{H}^{-2}((\a^{-2}_H(n_{(1)})1_3 \a_H^{-1}(m_{(1)}))S^{-1}(1_1))\\
&=& \a^{-1}_{N}(n_{(0)(0)})\o \a^{-2}_H(n_{(0)(1)})(1_2\c \a_M^{-2}(m_{(0)}))\o \a_{H}^{-2}(n_{(1)}[1_3 \a_H^{-2}(m_{(1)})]S^{-1}(1_1))\\
&=& \a^{-1}_{N}(n_{(0)(0)})\o \a_H^{-2}(n_{(0)(1)})\c \a^{-1}_M(m_{(0)})\o \a_{H}^{-2}(n_{(1)}m_{(1)})\\
&=&(c_{M,N}\o id_H)\rho_{M\o N}(m\o n).
\end{eqnarray*}
As for (4.2), for any $m\in M, n\in N$ and $p\in P$,  we have
\begin{eqnarray*}
&&(a_{N,P,M}\circ c_{M,N \o P}\circ a_{M,N,P})((m \o_t n) \o_t p)\\
&=&a_{N,P,M}(\a_{N}^{-1}(n_{(0)})\o p_{(0)}\o \a^{-4}_H(\a_H(p_{(1)})n_{(1)})\c \a^{-2}_{M}(m))\\
&=&\a_{N}^{-2}(n_{(0)})\o (p_{(0)}\o \a^{-3}_H(\a_H(p_{(1)})n_{(1)})\c \a^{-1}_{M}(m))\\
&=&(id_N \o_t (c_{M,P}))(\a_{N}^{-2}(n_{(0)})\o ( \a^{-2}_H(n_{(1)})\c\a^{-1}_{M}(m) \o \a_P(p))\\
&=&((id_N \o_t c_{M,P})\circ a_{N,M,P}\circ (c_{M,N} \o_t id_P))((m \o_t n) \o_t p),
\end{eqnarray*}
we can check that  (4.1) in the similar way.

{\bf Lemma 4.3.} $c_{M,N}$ is bijective with inverse
$$c_{M,N}^{-1}(n\o m)=\a_M^{-1}(\a_H^{-1}(S(n_{(1)}))\c m)\o \a_{N}^{-1}(n_{(0)}).$$

{\it Proof.} First, we prove that $c_{M,N}\circ c_{M, N}^{-1}=id$. For any $m\in M$ and $n\in N$,  we have
 \begin{eqnarray*}
    && c_{M,N}c_{M,N}^{-1}(n \o m)\\
    &=& c_{M,N}(\a_M^{-1}(\a_H^{-1}(S(n_{(1)}))\c m)\o \a_{N}^{-1}(n_{(0)}))\\
      &=& c_{M,N}(\a_H^{-2}(S(n_{(1)}))\c \a_M^{-1}(m)\o \a_{N}^{-1}(n_{(0)}))\\
      &=& \a_{N}^{-2}(n_{(0)(0)})\o \a_{N}^{-3}(n_{(0)(1)})\c [\a_H^{-3}(S(n_{(1)}))\c \a_M^{-2}(m)]\\
      &=& \a_{N}^{-2}(n_{(0)(0)})\o [\a_{N}^{-4}(n_{(0)(1)})\a_H^{-3}(S(n_{(1)}))]\c \a_M^{-1}(m)\\
   &=& \a_{N}^{-1}(n_{(0)})\o [\a_{N}^{-4}(n_{(1)1})\a_H^{-4}(S(n_{(1)2}))]\c \a_M^{-1}(m)\\
   &=&  \a_{N}^{-1}(n_{(0)})\o 1_2\varepsilon(1_1n_{(1)})\c \a_M^{-1}(m)\\
    &=&  1_1'\a_{N}^{-2}(n_{(0)})\o 1_2\varepsilon(1_11_2'n_{(1)})\c \a_M^{-1}(m)\\
    &=& 1_1\a_{N}^{-1}(n)\o 1_2\c \a_M^{-1}(m)\\
    &=&n\o m.
\end{eqnarray*}
Then, we note that the following relation holds, for all $m\in M$,
\begin{eqnarray*}
&&m_{(0)}\o m_{(1)}\\
   &=&\a_H^{-1}(1_{21})\c \a_M^{-1}(m_{(0)})\o[\a_H^{-2}(1_{22})\a_H^{-2}(m_{(1)})]S^{-1}(1_{1})\\
   &=& 1_{1}'\c (1_{2}\c \a_M^{-1}(m_{(0)}))\o[\a_H^{-2}(1'_{2})\a_H^{-2}(m_{(1)})]S^{-1}(1_{1})\\
   &=& 1_{1}'\c (1_{2}\c \a_M^{-1}(m_{(0)}))\o\a_H^{-1}(1'_{2})[\a_H^{-2}(m_{(1)})S^{-1}(\a_H^{-1}(1_{1}))]\\
   &=& 1_{1}'\c (1_{2}\c \a_M^{-1}(m_{(0)}))\o 1'_{2}[\a_H^{-2}(m_{(1)})S^{-1}(1_{1})]\\
   &=&1_{2}\c \a_M^{-1}(m_{(0)})\o \a_H^{-2}(m_{(1)})S^{-1}(1_{1}).
   \end{eqnarray*}
 Finally, we check that  $c_{M,N}^{-1}\circ c_{M, N}=id$. For any $m\in M$ and $n\in N$,  we have
\begin{eqnarray*}
&& c_{M,N}^{-1}c_{M,N}(m\o n)\\
 &=&c_{M,N}^{-1}(\a^{-1}_{N}(n_{(0)})\o \a_M^{-1}(\a^{-1}_H(n_{(1)})\c m))\\
&=& \a^{-3}_{H}(S(n_{(0)(1)}))\c [\a^{-3}_H(n_{(1)})\c \a_M^{-2}(m)]\o \a^{-2}_{N}(n_{(0)(0)})\\
&=& [\a^{-4}_{H}(S(n_{(0)(1)}))\a^{-3}_H(n_{(1)})]\c \a_M^{-1}(m)\o \a^{-2}_{N}(n_{(0)(0)})\\
&=& [\a^{-4}_{H}(S(n_{(1)1}))\a^{-4}_H(n_{(1)2})]\c \a_M^{-1}(m)\o \a^{-1}_{N}(n_{(0)})\\
&=& \varepsilon_s(n_{(1)})\c  \a_M^{-1}(m)\o \a^{-1}_{N}(n_{(0)})\\
&=& 1_1S ((S^{-1}\varepsilon_s(n_{(1)})))\c  \a_M^{-1}(m)\o 1_2 \c \a^{-1}_{N}(n_{(0)})\\
&=& S(1_2)S ((S^{-1}\varepsilon_s(n_{(1)}))) \c  \a_M^{-1}(m)\o S(1_1)\c \a^{-1}_{N}(n_{(0)})\\
&=& S ((S^{-1}\varepsilon_s(n_{(1)}))1_2) \c  \a_M^{-1}(m)\o S(1_1)\c \a^{-1}_{N}(n_{(0)})\\
&=& S ((S^{-1}\varepsilon_s(n_{(1)}))S(1_1)) \c  \a_M^{-1}(m)\o S^{2}(1_2)\c \a^{-1}_{N}(n_{(0)})\\
&=& S^{2}(1_1) \c  \a_M^{-1}(m) \o S^{2}(1_2S^{-1}\varepsilon_s(n_{(1)}))\c \a^{-1}_{N}(n_{(0)})\\
&=& S^{2}(1_1) \c  \a_M^{-1}(m) \o S^{2}(1_2)\c \a^{-1}_{N}(n)\\
&=& m\o n.
 \end{eqnarray*}

{\bf Theorem 4.4.} ${}_H{\cal WYD}^H$ is a braided monoidal category.

We can make now the connection between  Yetter-Drinfeld modules over weak Hom-Hopf algebras  and modules over quasitriangular
weak Hom-Hopf algebras.

{\bf Definition 4.5.}$^{\cite{ZW}}$  Let $(H,\a)$ be a weak Hom-bialgebra. If there exists $R = R^{(1)} \o R^{(2)} \in \Delta^{op}(1)(H \o_\Bbbk H)\Delta(1)$, such that the following conditions hold:\\
$\left\{\begin{array}{l}
(1)~~(\a \o \a)R =R;\\
(2)~~R\Delta(h) = \Delta^{op}(h)R;\\
(3)~~\mbox{there~exists~}\overline{R} \in \Delta(1)(H \o_\Bbbk H)\Delta^{op}(1),~\mbox{such~that~}R\overline{R} = \Delta^{op}(1),~\overline{R}R = \Delta(1);\\
(4)~~\a(R^{(1)}) \o R^{(2)}_1 \o R^{(2)}_2 = \a^{-1}(r^{(1)}R^{(1)}) \o R^{(2)} \o r^{(2)};\\
(5)~~R^{(1)}_1 \o R^{(1)}_2 \o \a(R^{(2)}) = r^{(1)} \o R^{(1)} \o \a^{-1}(r^{(2)}R^{(2)}),
\end{array}\right.$
\\
where $h \in H$, $r=R = R^{(1)} \o R^{(2)} = r^{(1)} \o r^{(2)}$, then $R$ is called an \emph{$R$-matrix} of $H$, $\overline{R}$ is called the \emph{weak inverse of $R$}. $(H,R)$ is called a \emph{quasitriangular weak Hom-bialgebra}.

{\bf Proposition 4.6.} Let $(H, \a_H, R)$ be a quasitriangular weak Hom-Hopf algebra, then we have

(i) Let $(M, \a_M)$ be a left $H$-module with action $H\o M\rightarrow M, h\o_t m\mapsto h\c m$. Define the linear map
$\rho_M: M\rightarrow M\o_t H, \rho_M(m) = m_{(0)} \o_t m_{(1)}:= R^{(2)} \c m\o_t \a_H(R^{(1)} )$. Then $(M, \a_M)$ with these structures
is a  Yetter-Drinfeld module over $H$.

(ii)  Let $(N, \a_N)$ be another left $H$-module with action $H\o_t N\rightarrow N, h\o_t n\mapsto h\c n$, regarded as a  Yetter-Drinfeld module as in (i), via the map $\rho_N: N\rightarrow N\o H, \rho_N(n)=n_{(0)} \o n_{(1)}:= r^{(2)}  \c m\o_t \a_H(r^{(1)} )$. We regard $(M \o_t N, \a_M\o\a_N)$ as a left $H$-module via the standard
action $h\c (m\o n)=h_1\c m\o_t h_2\c n$ and then we regard $(M \o_t N, \a_M\o \a_N)$ as a  Yetter-Drinfeld module
as in (i). Then this  Yetter-Drinfeld module $(M \o_t N, \a_M\o \a_N)$ coincides with the Yetter-Drinfeld
module $M\o_t N$ defined as in Proposition 4.1.

{\it Proof.} First we have to prove that $(M, \a_M)$ is a right $H$-comodule; $\rho(\a(m)) = \a_M(m_{(0)})\o \a(m_{(1)})$ is easy and left to the
reader, we check
\begin{eqnarray*}
(\a_M\o \Delta)\rho_M(m)&=& \a_M(R^{(2)}  \c m)\o \Delta(\a_H(R^{(1)} ))\\
&=& \a_H(R^{(2)} ) \c \a_M(m)\o \a_H(R^{(1)} _1)\o \a_H(R^{(1)} _2)\\
&=& r^{(2)}R^{(2)}  \c \a_M(m)\o \a^2_H(r^{(1)} )\o \a^2_H(R^{(1)} )\\
&=& \a_H(r^{(2)})\c (R^{(2)}  \c m)\o \a^2_H(r^{(1)} )\o \a^2_H(R^{(1)} )\\
&=& r^{(2)}\c (R^{(2)}  \c m)\o \a_H(r^{(1)} )\o \a^2_H(R^{(1)} )\\
&=& \rho_M(R^{(2)}  \c m)\o \a^2_H(R^{(1)} )\\
&=&(\rho_M\o \a_M)\rho_M(m).
\end{eqnarray*}

Now we check the  Yetter-Drinfeld condition (3.3):
\begin{eqnarray*}
(h_{2}\c m)_{(0)} \o (h_{2}\c m)_{(1)}\a_H^{2}(h_{1})&=& R^{(2)}\c (h_2\c  m)\o \a_H(R^{(1)})\a_H^{2}(h_{1})\\
&=& \a_H(R^{(2)})\c (h_2\c  m)\o \a^2_H(R^{(1)})\a_H^{2}(h_{1})\\
&=& (R^{(2)}h_2)\c  \a_M( m)\o \a^2_H(R^{(1)}h_{1})\\
&=& (h_1R^{(2)})\c  \a_M( m)\o \a^2_H(h_2R^{(1)})\\
&=& \a_H(h_1)\c (R^{(2)}\c m)\o \a^2_H(h_2) \a^2_H(R^{(1)})\\
&=&\a_H(h_{1})\c m_{(0)}\o \a_H^{2}(h_{2})\a_H(m_{(1)}).
\end{eqnarray*}

(ii) We only need to prove that the two comodule structures on $M\o_t N$ coincide, that is, for all $m\in M$ and $n\in N$,

\begin{eqnarray*}
m_{(0)} \o n_{(0)} \o \a^{-2}_H(n_{(1)}m_{(1)})=R^{(2)}\c (m\o n)\o \a_H(R^{(1)}),
\end{eqnarray*}
that is
\begin{eqnarray*}
R^{(2)} \c m \o r^{(2)}\c n \o \a^{-2}_H(\a_H(r^{(2)})\a_H(R^{(2)}))=R^{(2)}_1\c m \o R^{(2)}_2\c n \o \a_H(R^{(1)}),
\end{eqnarray*}
which is equivalent to
\begin{eqnarray*}
\a_H(R^{(2)}) \c m \o \a_H(r^{(2)})\c n \o r^{(2)}R^{(2)}=R^{(2)}_1\c m \o R^{(2)}_2\c n \o \a_H(R^{(1)}).
\end{eqnarray*}

{\bf Proposition 4.7.}    Let $(H, \a_H, R)$ be a quasitriangular weak Hom-Hopf algebra. Denote by $Rep(H)$ the category
whose objects are left $H$-modules and whose morphisms are $H$-linear maps.  Then $Rep(H)$ is a braided
monoidal subcategory of ${}_H{\cal WYD}^H$, with tensor product defined as in Proposition 4.1,
associativity constraints defined by the formula    $a_{M,N,P}((m \tilde{\o} n) \tilde{\o} p) = \a^{-1}_M(m) \tilde{\o}(n \tilde{\o} \a_P(p)) $
for any $M,N,P \in Rep(H)$, and braiding $c_{M,N}:M \o_t M\rightarrow N \o_t M,~~m \o_t n\mapsto R^{(2)}\cdot \a_N^{-1}(n) \o_t R^{(1)}\cdot \a_M^{-1}(m),$
with inverse
$c^{-1}_{M,N}: N \o_t M\rightarrow M \o_t M,~~n \o_t m\mapsto \overline{R}^{(1)}\cdot \a_M^{-1}(m) \o_t \overline{R}^{(2)}\cdot \a_N^{-1}(n),$
for any $(M,\a_M),(N,\a_N) \in Rep(H)$.

\section*{5. A braided monoidal category ${}_H{\cal WYD}^H$ II}
\def\theequation{5. \arabic{equation}}
\setcounter{equation} {0} \hskip\parindent

Modules over quasitriangular weak Hom-Hopf algebras become  Yetter-Drinfeld modules over weak Hom-Hopf algebras are proved in Section 4. Similarly, comodules over coquasitriangular weak Hom-Hopf algebras become Yetter-Drinfeld modules over weak Hom-Hopf algebras; inspired by this, we can introduce a second braided monoidal category
structure on ${}_H{\cal WYD}^H$. We include these facts here for completeness. Each of the next results is the analogue of a result in Section 4; their proofs are similar to those of their analogues and are left to the
reader.

{\bf Proposition 5.1.} Let $(H,\a_H)$ be a weak Hom-Hopf algebra.

(i) Let $(M,\a_M),  (N,\a_N)  \in
{}_H{\cal WYD}^H$, with notation as above,  and the tensor product $M \tilde{\o} N$ is obtained by
$$
M \tilde{\o} N = \{m \tilde{\o} n = m_0 \o_\Bbbk n_0 \varepsilon(m_1n_1)\mid m \in M, n \in N\},
$$ with structures:
\begin{eqnarray*}
&&h\c (m \tilde{\o} n)=\a_H^{-2}(h_{1})\c m \tilde{\o}\a_H^{-2}(h_{2})\c n,\\
&&m\tilde{\o} n \mapsto(m\tilde{\o} n)_{(0)}\o (m\tilde{\o} n)_{(1)}=(m_{(0)}\tilde{\o} n_{(0)})\tilde{\o} n_{(1)}m_{(1)}.
\end{eqnarray*}

(ii) ${}_H{\cal WYD}^H$ is a braided monoidal category, with tensor product $\tilde{\o}$ as in (i) and associativity
constraints $a_{M, N, P}$ and quasi-braiding $c_{M,N}$ defined as follows:
for any $(M, \a_M)$, $(N,\a_N), (P,\a_P) \in {}_H{\cal WYD}^H$, define an associativity constraint by
$$
a_{M,N,P}((m \tilde{\o} n) \tilde{\o} p) = \a^{-1}_M(m) \tilde{\o}(n \tilde{\o} \a_P(p)),~m \in M,~n \in N,~p \in P,
$$
$$ c_{M,N}: M \tilde{\o} N \rightarrow M \tilde{\o} M,\,\,\, c_{M,N}(m\tilde{\o}
n)=\a^{-1}_{N}(n_{(0)})\tilde{\o} \a_M^{-1}(\a^{-1}_H(n_{(1)})\c m).$$ with inverse
$$c_{M,N}^{-1}(n\tilde{\o} m)=\a_M^{-1}(\a_H^{-1}(S(n_{(1)}))\c m)\tilde{\o} \a_{N}^{-1}(n_{(0)})$$.

{\bf Definition 5.2.}$^{\cite{ZW}}$  Let $(H,\a)$ be a weak Hom-bialgebra, if there is a linear map $\sigma:H \o_\Bbbk H\rightarrow \Bbbk$, such that the following conditions hold:\\
$\left\{\begin{array}{l}
(1)~~\sigma(a,b) = \varepsilon(b_1a_1)\sigma(a_2,b_2)\varepsilon(a_3b_3);\\
(2)~~\sigma(a_1,b_1)a_2b_2 = b_1a_1\sigma(a_2,b_2);\\
(3)~~\mbox{there~exists~}\sigma' \in hom_\Bbbk(H \o H,\Bbbk),~\mbox{such~that~}\sigma(a_1,b_1)\sigma'(a_2,b_2) = \varepsilon(ab),\\
~~~~~~\mbox{and~}\sigma'(a_1,b_1)\sigma(a_2,b_2) = \varepsilon(b a);\\
(4)~~\sigma(\a(a),\a(b)) = \sigma(a,b);\\
(5)~~\sigma(\a(a),bc) = \sigma(a_1,\a(c))\sigma(a_2,\a(b));\\
(6)~~\sigma(ab,\a(c)) = \sigma(\a(a),c_1)\sigma(\a(b),c_2),
\end{array}\right.$
\\
where $a,b,c \in H$, then $\sigma$ is called an \emph{coquasitriangular form} of $H$, $\sigma'$ is called the \emph{weak convolution inverse of $\sigma$}. $(H,\sigma)$ is called a coquasitriangular weak Hom-bialgebra.

{\bf Proposition 5.3.} Let $(H, \a_H, \sigma)$ be a coquasitriangular weak Hom-Hopf algebra, then we have

(i) Let $(M, \a_M)$ be a left $H$-comodule with coaction $M\rightarrow M\tilde{\o} H, \rho_M(m) = m_{(0)} \tilde{\o} m_{(1)}$. Define the linear map
$H\o M\rightarrow M, h\tilde{\o} m \mapsto h\c m:= \sigma(\a_H(h),m_{(1)})m_{(0)}$. Then $(M, \a_M)$ with these structures
is a  Yetter-Drinfeld module over $H$.

(ii)  Let $(N, \a_N)$ be another left $H$-comodule with coaction $N\rightarrow N\tilde{\o} H, \rho_N(n) = n_{(0)} \tilde{\o} n_{(1)}$. Define the linear map
$H\o N\rightarrow N, h\tilde{\o} n \mapsto h\c n:= \sigma(\a_H(h),n_{(1)})n_{(0)}$. We regard $(M \tilde{\o} N, \a_M\o \a_N)$ as a left $H$-module via the standard
action $h\c (m\tilde{\o} n)=\a_{H}^{-2}(h_1)\c m\tilde{\o} \a_{H}^{-2}(h_2)\c n$ and then we regard $(M \tilde{\o} N, \a_M\o \a_N)$ as a  Yetter-Drinfeld module
as in (i). Then this Yetter-Drinfeld module $(M \tilde{\o} N, \a_M\o \a_N)$ coincides with the Yetter-Drinfeld
module $M\tilde{\o} N$ defined as in Proposition 5.1.

{\bf Theorem 5.4.} Let $(H, \a_H, \sigma)$ be a coquasitriangular weak Hom-Hopf algebra. Denote by $Corep(H)$ the category whose objects are right $H$-comodules $(M,
\a_M)$ and morphisms are morphisms of right $H$-comodules.  Then $Corep(H)$ is a braided
monoidal subcategory of ${}_H{\cal WYD}^H$, with tensor product defined as in Proposition 5.1,
associativity constraints defined by the formula    $a_{M,N,P}((m \tilde{\o} n) \tilde{\o} p) = \a^{-1}_M(m) \tilde{\o}(n \tilde{\o} \a_P(p)) $
for any $M,N,P \in Corep(H)$, and braiding $c_{M,N}: M \tilde{\o} N\rightarrow N\tilde{\o} M,~~m \tilde{\o} n\mapsto
\a_M^{-1}(n_0) \tilde{\o}\a_M^{-1}(m_0)\sigma(m_1, n_1),$
with inverse
$c^-1{}_{M,N}:N\tilde{\o} M\rightarrow M\tilde{\o}N,~~n \tilde{\o} m\mapsto\a_M^{-1}(m_0) \tilde{\o} \a_N^{-1}(n_0)\sigma'(m_1, n_1),$
for any $(M,\a_M),(N,\a_N) \in Corep(H)$.
\begin{center}
 {\bf ACKNOWLEDGEMENT}
\end{center}

   The work is supported by Project Funded by China Postdoctoral Science Foundation (No. 2015M580508),  the Fund of Science and Technology Department of Guizhou Province (No. 2014GZ81365) and the Program for Science and Technology Innovation Talents in Education Department of Guizhou Province(No. KY[2015]481).

\end{document}